\documentclass[a4paper]{article}
\usepackage{amsthm}
\usepackage{amsmath}
\usepackage{amssymb}
\usepackage{ifthen}
\usepackage{supertabular}
\newtheorem{theorem}{Theorem}[section]
\newtheorem{proposition}[theorem]{Proposition}
\newtheorem{lemma}[theorem]{Lemma}
\newtheorem{corollary}[theorem]{Corollary}
\newtheorem{remark}[theorem]{Remark}
\newtheorem{example}[theorem]{Example}

\begin{document}
\title{Structure theorems for AP rings}
\author{David W. Lewis - Stefan A.G. De Wannemacker}
\date{January, 2007}
\maketitle
\begin{abstract}In ``New Proofs of the structure
theorems for Witt Rings'', the first author shows how the standard
ring-theoretic results on the Witt ring can be deduced in a quick
and elementary way from the fact that the Witt ring of a field is
integral and from the specific nature of the explicit annihilating
polynomials he provides. We will show in the present article that
the same structure results hold for larger classes of commutative
rings and not only for Witt rings. We will construct annihilating
polynomials for these rings.
\end{abstract}
%%%%%%%%%%%%%%%%%%%%%%%%%%%%%%%%%%%%%%%%%%%%%%%%%%%%%%%%%%%%%%%%%%%%%%%%%%%%%%%%%%%%%%%%%%%%%%%%%%%%%%%%%%%%%%%%%%
%                                                                                                                %
%  SECTION : INTRODUCTION                                                                                        %
%                                                                                                                %
%%%%%%%%%%%%%%%%%%%%%%%%%%%%%%%%%%%%%%%%%%%%%%%%%%%%%%%%%%%%%%%%%%%%%%%%%%%%%%%%%%%%%%%%%%%%%%%%%%%%%%%%%%%%%%%%%%
\section{Introduction}
Using his theory of multiplicative forms, Pfister obtained in 1966
several structure theorems for Witt rings of quadratic forms over
fields. The proofs of these results were later simplified and
Harrison, Leicht and Lorenz added different results concerning the
ideal theory of these rings.  In 1987, the first author produced
polynomials annihilating all non-singular quadratic forms of a given
dimension in the Witt ring. He showed that the structure results
also follow from the fact that the Witt ring is integral and the
specific form of the polynomials. The main goal of this article is
to obtain these structure theorems for more general classes of
rings.

In the second section we will look at commutative rings $R$
additively generated by a subset of this ring $R$. The elements of
this subset being zeros of a polynomial $q(X)$ with integer
coefficients. In these rings $R$, which we will call annihilating
polynomial rings or AP rings for short, one can define a length map
$\ell:\,R\rightarrow\mathbb{N}$ and construct polynomials $p_n(X)\in
\mathbb{Z}[X]$ in such a way that all elements of length $n$ are
annihilated by $p_n(X)$. We will give many examples of AP rings.

In the third section we will obtain structure results for all AP
rings and in the fourth section we will look at AP rings with
generating polynomial $q(X)=X^{2^n}-1$. We will show a close
relationship between the signatures and the prime ideals in these
rings and assuming an extra 'admissible' condition on AP rings we
will be able to produce the complete classification of the spectrum,
in analogy to the work of Harrison, Leicht and Lorenz on the
spectrum of Witt rings. Finally, looking at the specific nature of
the prime ideals in these AP rings, we will obtain an analogue of
Pfister's local-global principle.

In the last section we will construct the annihilating polynomials
for different choices of $q(X)$.\\
\newpage
\section{Annihilating polynomial rings}
\subsection{Definition}
Let $q(X)$ be a monic polynomial with integer coefficients.  Let $R$
be a commutative ring and $S$ a subset, such that
\begin{enumerate}
    \item[(R1)] $\ R$ is additively generated by the subset $S\subseteq R$ and\\
    \item[(R2)] $\ q(s)=0$ for all $s \in S$.\\
\end{enumerate}
A commutative ring $R$ satisfying the conditions $(R1)$ and $(R2)$
will be called an \emph{ annihilating polynomial ring} or \emph{AP
ring} for short. $S$ will be called a \emph{generating set} and
$q(X)$ a \emph{generating polynomial} for the AP ring $R$.\\
\ \\
The condition $(R1)$ enables us to introduce a notion of length.
Given an element $r\in R$ we denote by $\ell_{S}(r)$ the least
number $l$ such that it is possible to write $r$ as a sum of $l$
elements of $S$, i.e. $r=\sum_{i=1}^l \epsilon_i a_i$ with
$\epsilon_i = \pm 1, a_i \in S$; we call $\ell_{S}(r)$ the
\emph{length of $r$ (with respect to $S$)}. In this way, we have
defined a map $\ell_{S} : R\longrightarrow
\mathbb N$ which we shall call the \emph{length map} of $R$ relative to $S$.\\
We will write $\ell(r)$ for short if the choice of the generating
set $S$ for the ring $R$ is clear.\\
\ \\
Recall this general proposition for commutative
rings. \\
Let $q_1(X),q_2(X),\ldots,q_n(X)$ be monic polynomials with integer
coefficients and such that none of the $q_i(X)$ has a multiple root
in the complex numbers $\mathbb{C}$.  Write $R_i$ for the set of
roots of $q_i(X)$ in $\mathbb{C}$ and $T_n$ for the set of all
complex numbers expressible in the form $\sigma =
\sum_{i=1}^n\epsilon_i\sigma_i$ where $\epsilon_i=\pm 1,\
\sigma_i\in R_i$ for each $i$.  We write $p_n(X)=\prod_{\sigma\in
T_n}(X-\sigma)$ which is a monic polynomial with integer
coefficients, without multiple roots. We have the following
\begin{proposition}\label{prop_lewis}[Lewis]
Let $R$ be a commutative ring.  Let $x_1,x_2,\ldots,x_n$ be elements
of $R$ such that $q_i(x_i)=0$ for $i=1,2,\ldots,n$ and let
$x=\sum_{i=1}^n\epsilon_ix_i$ where $\epsilon_i=\pm 1$.  Then
$p_n(x)=0$ in $R$.
\end{proposition}
\begin{proof}
See \cite{LEW90}.
\end{proof}
Let $R$ be an AP ring with generating polynomial $q(X)$. Let $T_1$
denote the set of roots of $q(X)$ in $\mathbb{C}$ and $T_n$ the set
of complex numbers expressible in the form
$\sigma=\sum_{i=1}^n\epsilon_i\sigma_i$ where $\epsilon_i=\pm 1,\
\sigma_i\in T_1$. Put $p_n(X)= \prod_{\sigma\in T_n}(X-\sigma)$. We
have:
\begin{corollary}
Let $R$ be an AP ring. Every element of length $n$ is annihilated
by $p_n(X)$ in $R$.  In particular, $R$ is integral.
\end{corollary}
\begin{proof}
Let $r$ be an element of length $n$ in the AP ring $R$ with
generating set $S$ and generating polynomial $q(X)$. Then there
exist elements $a_1,a_2,\ldots,a_n\ \in S$ such that
$r=\sum_{i=1}^n\epsilon_ia_i,\ \epsilon_i=\pm1$ and
$q(a_1)=q(a_2)=\ldots=q(a_n)=0.$ From the previous proposition it
follows that $p_n(r)=0.$
\end{proof}
Before giving examples of AP rings, let us give an example of
constructing an annihilating polynomial.\\
Let $R$ be an AP ring with generating polynomial $q(X)=X^2 - 1.$ The
roots of the generating polynomial are $-1$ and $1$. The possible
values for the sum of $n$ elements out $\{-1,1\}$ lie in
$\{-n,-n+2,\ldots,n-2,n\}$. The annihilating polynomial for an
element of length $n$ is thus the $n$-th Lewis polynomial
\[p_n(X) = (X-n)(X-(n-2))\ldots(X+(n-2))(X+n).\]

Further examples of constructing the annihilating polynomial for a
given AP ring will be given in section 5.
\subsection{Examples of AP rings}
\begin{enumerate}
\item[(i)]$\mathbb{Z}$ is an AP ring.  It is additively generated
by $S=\{-1,1\}$ and we can take the generating polynomial
$q(X)=X^2-1$.\\
\item[(ii)] A product $\prod\mathbb{Z}$ of finitely many copies of
$\mathbb{Z}$ is an AP ring.\\  It is generated by the elements $S=\{
\pm e_i\}$ where $e_i$ has a $1$ in the $i$-th place and
zero elsewhere, and we can take $q(X)=X^3-X$ as a generating polynomial.\\
\item[(iii)]Let $G$ be an abelian group of finite exponent $n$.
The group ring $\mathbb{Z}[G]$ is an AP ring with generating set
$S=G$. Since every element has finite order $n$, we can take
$q(X)=X^{n}-1$ as the generating
polynomial.\\
\item[(iv)]Let $K$ be an ideal in the AP ring $R$. Then $R/K$ is a
AP ring as well, with generating set $S/K$.  If $q(X)$ is a
generating polynomial for the AP ring $R$, then $q(X)$ is a
generating polynomial for $R/K$ as well. Examples are
$\mathbb{Z}_n$, abstract Witt rings (see \cite{Marshall}),\ldots \\
\item[(v)] Witt-Grothendieck $\hat{W}(F)$ and Witt rings $W(F)$ of
quadratic forms over fields $F$ (see \cite{LAM} or \cite{SCH} for
further details on quadratic forms and Witt rings).\\ Let $F$ be a
field of characteristic not 2. Consider the group ring
$\mathbb{Z}[F^\ast/F^{\ast^2}]$ and the canonical ring homomorphism
\[\pi : \mathbb{Z}[F^\ast/F^{\ast^2}] \longrightarrow \hat{W}(F)\]
defined by $\alpha \mapsto \langle \alpha\rangle$ for $\alpha \in
F^\ast/F^{\ast^2}$. $\pi$ is surjective since every bilinear space
is an orthogonal sum of $1$-dimensional ones and therefore
\[\hat{W}(F)\cong\mathbb{Z}[F^\ast/F^{\ast^2}]/{\rm ker}(\pi).\]
The fact that $\hat{W}(F)$ is an AP ring follows from examples
$(iii)$ and $(iv)$ and we can take $q(X)=X^2-1$.\\
$W(F)$ is the quotient ring of $\hat{W}(F)$ by the ideal generated
by the hyperbolic spaces $\{ n[\mathbb{H}]\,\mid\,n \in
\mathbb{Z}\}$. By example $(iv)$, $W(F)$ is an AP ring with
generating polynomial $q(X)=X^2-1$.
\item[(vi)] Witt rings of higher level (see \cite{Kleinstein Rosenberg}).\\
Let $F$ be a field of characteristic not $2$. Denote by $\langle
\alpha \rangle_n$ the class $\alpha F^{\ast^{2^n}}$ in
$F^{\ast}/F^{\ast^{2^n}}$ and consider the ideal $R_n(F)$ of
$\mathbb{Z}[F^{\ast}/F^{\ast^{2^n}}]$ generated by $\langle 1
\rangle_n + \langle -1\rangle_n$ and $(\sum_{i=0}^{2^n-1}\langle
\alpha^i \rangle_n)(1-\langle1+\alpha\rangle_n)$ for all $\alpha \in
F^{\ast}$.\\ The ring
$W_n(F):=\mathbb{Z}[F^{\ast}/F^{\ast^{2^n}}]/R_n(F)$ is called the
{\em Witt ring of level $n$ of $F$} and is an AP ring with
generating polynomial $q(X)=X^{2^{n}}-1$.
\item[(vii)] Products and tensor products of finitely many AP rings are
AP rings.  If $q_1(X), q_2(X), \ldots, q_n(X)$ are generating
polynomials for respectively\\ $R_1, R_2, \ldots, R_n$ then $q(X) =
\prod_{i=1}^nq_i(X)$ is a generating polynomial for $R_1 \times R_2
\times \ldots\times R_n$ and for $R_1 \otimes R_2
\otimes \ldots \otimes R_n$.\\
\item[(viii)] The Burnside ring of a finite group (see
\cite{TtD} for further details on Burnside rings). Let $G$ be a
finite group. The set of isomorphism classes of finite $G$-sets form
a commutative associative semi-ring $\Omega^{+}(G)$ with unit under
disjoint union and cartesian product.  The Grothendieck ring
$\Omega(G)$ constructed
from this semi-ring is called the {\em Burnside ring} of $G$.\\
Additively, $\Omega(G)$ is the free abelian group on isomorphism
classes of transitive $G$-sets. Equivalently, an additive
$\mathbb{Z}$-basis for $\Omega(G)$ is given by the $[G/H]$ in
$\Omega(G)$ where $(H)$ runs through the set $C(G)$ of conjugacy
classes of subgroups of $G$.\\ There exist an injective ring
homomorphism
\[ \varphi\ :\ \Omega(G) \longrightarrow \prod_{(H) \in C(G)}\mathbb{Z}\]

induced by
\[ T \longmapsto (|T^H|\ \mid\ (H) \in C(G))
\]
where $|T^H|$ denotes the number of elements of the $G$-set $T$ fixed under $H$.\\
Each generator $y=[G/H]$ of $\Omega(G)$ maps to an element
$(n_1,n_2,\ldots,n_k)\ $ ($k=|C(G)|$) and this is annihilated in
$\prod_{i=1}^k\mathbb{Z}$ by the polynomial
$q_y(X)=\prod_{i=1}^k(X-n_i)$.  By the injectivity of the above ring
homomorphism $\varphi$, $q_y(X)$ annihilates $y$ in $\Omega(G)$. Let
$q(X)=\prod q_y(X)$, the product taken over $S$, the finite set of
generators $y$ of $\Omega(G)$. Then $\Omega(G)$ is
an AP ring with generating polynomial $q(X)$.\\
\item[(ix)]Witt ring of a central simple algebra
In \cite{LEW_TIG} the first author and Tignol define this ring as a
quotient of a group ring of abelian groups of exponent two. They
generalize the notion of Witt ring of a skew field done earlier by
Craven and Sladek (see \cite{CRA} and \cite{SLA}). All these Witt
rings are AP rings with generating polynomial $q(X)=X^2-1$. Note
that this immediately implies that the Lewis polynomials annihilate
the elements of these rings.
\end{enumerate}
In the next section we will obtain structure results for these AP
rings using the fact that these rings are integral and the
specific nature of the annihilating polynomials.\\
\ \\
\newpage
\section{Structure theorems for AP rings}
We will start this section with investigating the spectrum of AP
rings.\\
Let us fix the following notations first.\\

${\rm Nil} R$ is the nilradical of $R$, i.e. set of nilpotent
elements of $R$, ${\rm Nil} R = \cap P$ where $P$ runs through all
the prime ideals of $R$.\\

$R_t$ is the torsion subgroup of the additive group of $R$.\\

${\rm Spec} R$ is the set of prime ideals of $R$.\\

${\rm Max} R$ and ${\rm Min} R$ are the subsets of ${\rm Spec} R$
consisting of the maximal ideals and minimal prime ideals of
$R$ respectively.\\

\begin{lemma}
Let $R$ be an AP ring. Then
\begin{enumerate}
\item[(i)] A prime ideal $P$ of  $R$ is maximal if and only if
there is a rational prime $p$ with $p\cdot 1_R \in P$. \item[(ii)]
If $R_t \subset {\rm Nil} R$ then a prime ideal $P$ of $R$ is
minimal if and only if $P \cap \mathbb{Z} = 0$, maximal otherwise,
and every maximal ideal of $R$ properly contains a minimal prime
ideal. \item[(iii)] ${\rm Nil} R \subset R_t$.
\end{enumerate}
\end{lemma}
\begin{proof}
\ \\
\begin{enumerate}
\item[(i),(ii)]Since every AP ring $R$ is integral over
$\mathbb{Z}$,  this follows from \cite[Lemma 2.5.]{KRW}.
\item[(iii)] Let $r \in {\rm Nil} R$ and let $k \in \mathbb{N}$ be
such that $r^k=0$ and $r^{k-1} \neq 0$. Since $R$ is an AP ring, $r$
is annihilated by some $p_n(X)$. Write
$r^n+a_{n-1}r^{n-1}+\ldots+a_1r+a_0 = 0$ where $a_i \in \mathbb{Z}$.
Suppose that $a_i$ is the non-zero number of lowest index ($i=0$ or
$i=1$ since $p_n(X)$ has no multiple roots (see \ref{prop_lewis})).
Multiplying the above equation with $r^{k-i-1}$ ($k-i-1 \geq 0$,
since $k\geq 2$) yields $a_ir^{k-1}=0$. Using this and multiplying
the equation with $a_ir^{k-i-2}$ we obtain $a_i^2r^{k-2}=0$.
Repeating this process will give $a_i^{k-1}r=0$, i.e. $r\in R_t$.
\end{enumerate}
\end{proof}
\ \\
Let $R$ be an AP ring with generating polynomial $q(X)$ and a multiplicatively closed generating set $S$.
Suppose $R$ is of the form $\bigoplus_{s\in S}\mathbb{Z}s$. Let $\chi$ be a (monoid-)morphism of $S$ into $\mathbb{C}$,
\[ \chi:\,S\longrightarrow \mathbb{C},\] mapping $1$ to $1$. Remark that
for every $s\in S$, $\chi(s)$ will be a root of $q(X)$. \\
Let $K$ be the field generated over $\mathbb{Q}$ by all the roots of
$q(X)$ (i.e. by all the $\chi(s)$) and let $\mathcal{C}$ be the
integral closure of $\mathbb{Z}$ in $K$.  Every morphism $\chi$ of
$S$ to $\mathbb{C}$ extends to a ring homomorphism
\[\phi_{\chi}:\,R\longrightarrow \mathcal{C}\] in the obvious
way by defining \[\phi_{\chi}(\sum a_is_i)=\sum
a_i\chi(s_i)\quad\text{for all}\quad a_i\in \mathbb{Z},\,s_i \in
S.\] The spectrum of $R=\bigoplus_{s\in S}\mathbb{Z}s$ is completely determined
by the following lemma.\\
\begin{lemma}
The minimal prime ideals of $R=\bigoplus_{s\in S}\mathbb{Z}s$ are the kernels
$P_{\chi}$ of the morphisms
$\phi_{\chi}:\,R\longrightarrow \mathcal{C}$. The
maximal ideals of $\mathbb{Z}[S]$ are of the form
$M_{\chi,\mathfrak{p}}=\phi_{\chi}^{-1}(\mathfrak{p})$ where
$\mathfrak{p}$ is a non-zero (and thus maximal) prime ideal of
$\mathcal{C}$.
\end{lemma}
\begin{proof}
The proof of \cite[Lemma 3.1.]{KRW} extends in an obvious way from
abelian torsion groups $G$ to sets $S$ and the assertion follows.
\end{proof}
\ \\
A \emph{signature} $\sigma$ of $R$ is a (necessarily surjective)
ring homomorphism $\sigma\ :\ R\longrightarrow \mathbb{Z}$. We
call a prime ideal $\mathcal{P}$ a \emph{signature ideal of }$R$
if $R/\mathcal{P}$ is isomorphic to $\mathbb{Z}$. We will call
$X(R)=\{\mathcal{P}\mid \mbox{$\mathcal{P}$ a signature ideal of
$R$}\}$ the \emph{space of signatures of }$R$.
\begin{proposition}
\ \\
Let $R=\bigoplus_{s\in S}\mathbb{Z}s$ be the AP ring generated by $q(X)$. If $q(X)$
has all its roots in $\mathbb{Z}$ then \[{\rm Spec}(R)={\rm Min}(R)
\cup {\rm Max}(R)\] where \[{\rm Min}(R)=X(R)\] and
\[{\rm Max}(R)=\{\mathcal{P} + pR\quad\mid\quad\mathcal{P}\in X(R),\quad\mbox{ and $p$
a rational prime}\}.\]
\end{proposition}
\begin{proof}
$q(X)$ has all its roots in $\mathbb{Z}$, so $\chi(s) \in
\mathbb{Z}$ for all $s \in S$.  Using the above notation we get that
$\mathcal{C}=\mathbb{Z}$ and $K=\mathbb{Q}$.  All the minimal prime
ideals, as kernels of homomorphisms of $R=\bigoplus_{s\in S}\mathbb{Z}s$
to $\mathbb{Z}$ are the signature ideals of $R$.\\
Since all the maximal ideals of $\mathbb{Z}$ are of the form
$p\mathbb{Z}$, the maximal ideals of $R$ are of the form
$\mathcal{P} + pR$, where $\mathcal{P}\in X(R).$
\end{proof}
The Burnside ring $\Omega(G)$ of a finite group $G$ satisfies the
conditions of the previous proposition (see example (viii)).  For
a subgroup $U\leq G$ and a $G$-set $S$ the map $S \mapsto
\#S^{U}$(i.e. the number of elements in $S$ invariant under $U$)
extends to a ring homomorphism $\phi_U : \Omega(G) \longrightarrow
\mathbb{Z}$.  Define for $p$ being 0 or a prime number
 the prime ideal $\mathfrak{p}_{_{U,p}}=\{ x \in \Omega(G)\,|\,\#x^{U}\equiv 0
\mod 2\}$. Since $X(\Omega(G))=\{{\rm ker}(\phi_U)\,\mid\,U\leq
G\}$ we can retrieve Dress' result on the description of the
spectrum of the Burnside ring $\Omega(G)$:\\
\begin{proposition}[Dress](see \cite[Proposition 1.]{Dress})
One has $\mathfrak{p}_{_{U,p}}\subset \mathfrak{p}_{_{V,q}}$ if and
only if $p=q$ and $\mathfrak{p}_{_{U,p}}= \mathfrak{p}_{_{V,q}}$ or
$p=0, q\neq 0$ and $\mathfrak{p}_{_{U,q}}\subset
\mathfrak{p}_{_{V,q}}.$ Especially, $\mathfrak{p}_{_{U,p}}$ is
minimal, respectively maximal, if and only if $p=0$, respectively
$p\neq 0.$
\end{proposition}
\newpage
\section{AP rings with generating polynomial $X^{2^k}-1$}
We will now take a closer look at AP rings with generating
polynomial $q(X)=X^{2^k}-1$ and a group $S=G$ as a set of
generators.  We will obtain some additional structure theorems.\\
Consider $I:=\langle 1-a\ |\ a\in G\rangle$.  We will call this
ideal the \emph{fundamental ideal of }$R$.\\
Observe that for any AP ring $R$, $\ |R/I|\leq 2$, since every
element is a sum of an even or an odd number of elements of $G \cup
-G$.
\begin{proposition}
Let $R$ be an AP ring with generating polynomial $X^{2^k}-1$ and a
group $S=G$ as a set of generators.  The following conditions are
equivalent :\\
\begin{enumerate}
\item[(i)] $\ R$ does not have odd characteristic. \item[(ii)] $\
|R/I|=2$. \item[(iii)] $\ \sum_{i=1}^m\varepsilon_i a_i = 0,\ a_i
\in G \Rightarrow\ m\in 2\mathbb{Z}$.
\end{enumerate}
\end{proposition}
\begin{proof}
\ \\
$(iii) \Rightarrow (ii) :$ Suppose $|R/I|=1$, then $1\in I$ and we
have $1+\sum_i\varepsilon_i(a_i-1) = 0$, contradicting $(iii)$.\\
$(ii) \Rightarrow (i) :$ Suppose $R$ has odd characteristic $k$
and let $r\in R$.  Then either $r\in I$ or $r=r+k \in I$ implying
$|R/I|=1$, a contradiction.\\
$(i) \Rightarrow (iii) :$ Assume $r=\sum_{i=1}^ma_i=0$ with $m$ odd.
Since $p_m(r)=0$ this implies that $p_m(0)=0$ in $R$ contradicting
the fact that $R$ has even characteristic. See also lemma
\ref{odd constant term}.\\
\end{proof}
\ \\
Any AP ring $R$ which satisfies the conditions of the previous
proposition will be called \emph{ admissible}.\\
\ \\These admissibility conditions are not always satisfied, e.g.
$\mathbb{Z}/n\mathbb{Z}$ is an AP ring for all $n\in \mathbb{N}$,
but is only admissible if
$n$ is even.\\
%\begin{remark}
\noindent Consider the Arason-Pfister property $AP(k)$
\[\text{ If}\quad  r = a_1 + a_2 + \ldots + a_n \in I^k,\ \text{with}\ a_i \in G\cup -G,\  n<2^k,\ \text{ then
}\ r=0.\] We have the following
\begin{proposition}
An AP ring is admissible if and only if $AP(1)$ holds.
\end{proposition}
\ \\
Let $R$ be an admissible AP ring.  The unique isomorphism
$R/I\cong \mathbb{Z}/2\mathbb{Z}$ induces a homomorphism
\begin{align*}
\wp : R&\rightarrow \mathbb{Z}/2\mathbb{Z} \\
 &r\mapsto \wp(r)\equiv \ell(r) \mod{2},\ \forall r \in R.
\end{align*}
\newpage
\begin{proposition}
Let $R$ be an admissible AP ring. Then the fundamental ideal $I$ is
the only prime ideal of index $2$ in $R$.
\end{proposition}
\begin{proof}
Suppose $\mathcal{P}$ is a prime ideal of index 2.  For every
generating element $a$ we have
$0=(1+a^{2^{k-1}})(1-a^{2^{k-1}})\in\mathcal{P}$. Since $2 \in
\mathcal{P}$, this implies that $1-a^{2^{k-1}} \in \mathcal{P}$.
Using this observation we obtain finally that $1-a \in \mathcal{P}$.
Since $1-a$ are exactly the generators for $I$, we get $I
\subseteq\mathcal{P}$ and thus $I=\mathcal{P}$ by maximality of
$I$.\end{proof}
\par
From now on we look at the special case $k=1$, i.e. where $R$ is an
AP ring with generating polynomial $q(X)=X^2-1$.
\begin{remark}
AP rings satisfying these conditions are for example Witt rings of
fields, Witt rings of central simple algebras and abstract Witt
rings.
\end{remark}
\begin{proposition}
Let $R$ be an AP ring with generating polynomial $q(X)=X^2-1$. If
$\mathcal{P}$ is a prime ideal of $R$, then either $\mathcal{P}$ is
of finite index which is a prime number, or $R/\mathcal{P}$ is
isomorphic to $\mathbb{Z}$
\end{proposition}
\begin{proof}
Let $\mathcal{P}$ be a prime ideal in $R$.  Since $R/\mathcal{P}$
is an integral domain and $(a-1)(a+1)=0$ for all $a \in G$, we
have that $a=1 + \mathcal{P}$ or $a=-1 + \mathcal{P}$ for all $a
\in G$.  Since $R$ is additively generated by $G$, we have that
$R/\mathcal{P}$ is cyclic generated by $1+\mathcal{P}$. If the
characteristic of the integral domain $R/\mathcal{P}$ is $p$, then
$\mathcal{P}$ is of finite index $p$ in $R$. Otherwise, the
characteristic of $R/\mathcal{P}$ is $0$ and so $R/\mathcal{P}$ is
isomorphic to $\mathbb{Z}$.
\end{proof}
\begin{proposition}
Let $R$ be an AP ring with generating polynomial $q(X)=X^2-1$.
There is a canonical one-to-one correspondence between signature
ideals in $R$ and the different signatures of $R$.
\end{proposition}
\begin{proof}
\ \\
Let $\sigma$ be a signature of $R$.  It is clear that ${\rm
ker}(\sigma)$ is a signature ideal in $R$.
\ \\
Conversely, suppose that $\mathcal{P}$ is a signature ideal of
$R$. The isomorphism $R/\mathcal{P}\cong \mathbb{Z}$ induces an
unique homomorphism $R\rightarrow\mathbb{Z}$ with kernel
$\mathcal{P}$.
\end{proof} \vspace*{2ex} \noindent
\ \\
Let $\mathcal{P}$ be a signature ideal of $R$ and consider the
obvious prime ideals of finite index $p$, namely $\mathcal{P}+pR$.
We will show that these ideals are the only prime ideals of finite
index in $R$.  Let us first recall the following result: \\
\begin{lemma}
Let $R$ be an AP ring with generating polynomial $q(X)=X^2-1$ and
assume $X(R)\neq \emptyset$. Then $\mathcal{P}+pR$ is the unique
prime ideal of finite index $p$, containing the signature ideal
$\mathcal{P}$.
\end{lemma}
\ \\
Moreover, these ideals are the only ideals of finite index $p$
different from 2 since:\\
\begin{lemma}
Let $R$ be an AP ring with generating polynomial $q(X)=X^2-1$ and
assume $X(R)\neq \emptyset$. Let $\mathcal{Q}$ be a prime ideal of
finite index $p\neq 2$. Then $\mathcal{Q} = \mathcal{P}+pR$ for some
signature ideal $\mathcal{P} \in X(R)$.
\end{lemma}
\begin{proof}
We will construct a signature $\sigma : R \longrightarrow
\mathbb{Z}$ such that $\mathcal{P}={\rm ker}(\sigma)$ and
$\mathcal{P}\subset\mathcal{Q}$.  The above lemma will then
complete the proof.\\
Let $\mathcal{M}$ be the minimal prime ideal such that
$\mathcal{M}\subset\mathcal{Q}$. Since $a^2 = 1$, for all $a \in G$,
we have $a-1 \in \mathcal{M}$ or $a+1 \in \mathcal{M}$. Define
$\sigma_{\mathcal{M}} : G\rightarrow \{-1,+1\}$ for all $a \in G$ by
\begin{equation}
\sigma_{\mathcal{M}}(a) =
\begin{cases} -1,&\text{if}\quad a+1 \in \mathcal{M}\\
+1,&\text{if}\quad a-1 \in \mathcal{M}.\end{cases}
\end{equation}
Define $\sigma : R \rightarrow \mathbb{Z}$ for all
$r=\sum_{i=1}^k\varepsilon_ia_i \in R$, where ${\varepsilon}_i=\pm
1,\ a_i \in G$ by
\[\sigma(r)=\sum_{i=1}^k\varepsilon_i\sigma_{\mathcal{M}}(a_i).\]
We claim that $\sigma$ is a well-defined signature of $R$. Proof
of claim:\\
Suppose $r = \sum_{i=1}^k\varepsilon_ia_i = 0 \in R$. Then
\begin{align*}
\sigma(r)&=\sum_{i=1}^k\varepsilon_i\sigma_{\mathcal{M}}(a_i)\\
&=\sum_{i=1}^k\varepsilon_i(a_i-\sigma_{\mathcal{M}}(a_i))\\
&\in \mathcal{M}\cap\mathbb{Z}. \\
\end{align*}
Since the AP ring $R$ is integral over $\mathbb{Z}$, any minimal
prime ideal of $R$ lies over $\mathbb{Z}$, i.e.
$\mathcal{M}\cap\mathbb{Z} = \{0\}.$ So $\sigma(r)=0$ and $\sigma$
is well-defined. Further, it is obvious that $\sigma(1)=1$ and
$\sigma(r_1+r_2)=\sigma(r_1)+\sigma(r_2)$, for all $r_1,r_2 \in R$,
from the definition of $\sigma$. To show that
$\sigma(r_1r_2)=\sigma(r_1)\sigma(r_2)$, for all $r_1,r_2 \in R$,
observe that it is sufficient to show that
$\sigma_\mathcal{M}(ab)=\sigma_\mathcal{M}(a)\sigma_\mathcal{M}(b)$,
for all $a,b \in G$.\\
$(a-\sigma_\mathcal{M}(a))(b-\sigma_\mathcal{M}(b))=ab-\sigma_\mathcal{M}(a)b-\sigma_\mathcal{M}(b)a+\sigma_\mathcal{M}(a)\sigma_\mathcal{M}(b)
\in \mathcal{M}$. Since
$\sigma_\mathcal{M}(a)b-\sigma_\mathcal{M}(a)\,\sigma_\mathcal{M}(b)
\in \mathcal{M}$ and
$\sigma_\mathcal{M}(b)a-\sigma_\mathcal{M}(b)\,\sigma_\mathcal{M}(a)
\in \mathcal{M}$, we have
$ab-\sigma_\mathcal{M}(a)\sigma_\mathcal{M}(b) \in \mathcal{M}$,
i.e. $\sigma_\mathcal{M}(a)\sigma_\mathcal{M}(b) =
\sigma_\mathcal{M}(ab)$.\\
Finally, consider the signature ideal $\mathcal{P}=\{r\in R | \sigma(r)=0\}$.\\
For $\,r=\sum_{i=1}^k\varepsilon_i a_i \in \mathcal{P}\,$, we have
$\,\sum_{i=1}^k\varepsilon_i\sigma_\mathcal{M}(a_i)=0\,$. This
implies (that $k$ is even and) that one half of the
$\,\varepsilon_i\sigma_\mathcal{M}(a_i)\,$ equals $-1$ and the other
half equals $+1$. Without loss of generality we can assume that
$\,\varepsilon_i\sigma_\mathcal{M}(a_i)=-1\,$ for $\,1\leq i \leq
k/2\,$ and $\,\varepsilon_i\sigma_\mathcal{M}(a_i)=+1\,$ for
$\,k/2 < i \leq k\,$. We can rewrite\\
\begin{align*}
r&=\sum_{i=1}^k\varepsilon_ia_i\\
&=\sum_{i=1}^{k}\varepsilon_i(a_i-\sigma_\mathcal{M}(a_i))\\
&\\
&\in \mathcal{M}\\
\end{align*}
So $\,\mathcal{P}\subset\mathcal{M}\subset\mathcal{Q}\,$.
\end{proof}
\ \\
This brings us to the complete classification of the prime ideal
spectrum of an admissible AP ring $R$, denoted ${\rm Spec}(R)$ :
\\
\begin{proposition}
\ \\
Let $R$ be an admissible AP ring with generating polynomial $q(X)=X^2-1$.\\
If $X(R)=\emptyset$, then ${\rm Spec} (R) = \{I\}$.\\
Otherwise ${\rm Spec} (R)\,=\, {\rm Min}(R)\,\cup\, {\rm Max} (R)$
where ${\rm Min} (R) = X(R)$ and ${\rm Max} (R) = \{I\}\cup
\{\mathcal{P} + pR\mid\mathcal{P}\in X(R), \mbox{$p$ odd
prime}\}$.\\
In particular, $R$ is a local ring if and only if
$X(R)=\emptyset$.
\ \\
\end{proposition}

The admissible AP ring $R$ is a local ring if, and only if, $I$ is
the only prime ideal in $R$. Otherwise a prime ideal of $R$ is
either a minimal prime ideal or a maximal ideal. The minimal prime
ideals are the ideals in ${\rm Min}(R)$. The maximal ideals are
the ones given by ${\rm Max}(R)$.

We will determine the following objects, which are very useful to
determine the structure of a ring $R$:\\
\begin{enumerate}
\item[$\star\ $] ${\rm Nil}(R)$, the nil radical of $R$,
consisting of all nilpotent elements of $R$, \item[$\star\ $] ${\rm
Jac}(R)$, the Jacobson radical of $R$, the intersection of all
maximal ideals of $R$, \item[$\star\ $] $R^\times$, the units (or
invertible elements) of $R$, \item[$\star\ $] ${\rm Zd}(R)$, the
zerodivisors in $R$ (including the zero element),
\item[$\star\ $] the idempotents in $R$ and \item[$\star\ $]
$R_t$, the torsion elements of $R$.
\end{enumerate}
\ \\
We will make a distinction between two cases, namely $X(R) =
\emptyset$ and $X(R)\neq\ \emptyset.$  From now on, $R$ will denote an admissible AP ring.\\
\subsection{$X(R) = \emptyset$}

In this case $I$ is the only prime (and thus the maximal) ideal in
$R$, which makes $R$ a local ring.\\
\ \\
Since ${\rm Nil}(R)$ is the intersection of all prime ideals in
$R$, we get ${\rm Nil}(R)=I$. Obviously, ${\rm Jac}(R)=I$, as it
is the only maximal ideal.\\
\ \\
In a local ring $R$, with unique maximal ideal $I$, the
multiplicative group of invertible elements consists of the
elements in  $R\setminus I$.\\
\ \\
For a commutative ring $R$, the set of zerodivisors, ${\rm Zd}(R)$,
is the union of a certain set of prime ideals in $R$. Given a
zerodivisor $z \in {\rm Zd}(R)$, the prime ideal $P$ containing $z$
is the prime ideal $P$,
maximal in the sense that $(R\setminus {\rm Zd}(R)) \cap P = \emptyset.$\\
That ${\rm Zd}(R)=I$ follows from the observation that $I={\rm
Nil}(R)\subseteq {\rm Zd}(R)$ and
that $R\setminus I = R^\times\subseteq R\setminus {\rm Zd}(R)$.\\
\ \\
The only idempotents are the trivial ones, i.e. 0 and $1$.  This
follows from the observation that for $e \in R$ one has $e \in I$
or $e-1 \in I$. Since ${\rm Zd}(R)=I$, we have $e(e-1)=0$ implying
$e-1=0$ or $e=0$ respectively.  So $R$ is a 'connected' ring.\\
\ \\
Since $I={\rm Nil}(R)$, we have $2 = 1 + 1 \in {\rm Nil}(R)$.  So
there exists a natural number $k$, such that $2^k=0$ in $R$.  This
implies that $2^k\,r=0,\ \forall r \in R$. So $R_t = R$.\\

These results can be summarized in the following
\begin{proposition}
Let $R$ be an admissible AP ring with generating polynomial $q(X)=X^2-1$ and assume that $X(R)=\emptyset$.\\
Then all elements in $R$ are torsion elements. Moreover, for an
element $r \in R$, the following conditions are equivalent :\\
\begin{itemize}
\item[(i)] $r\in I$, \item[(ii)] $r$ is nilpotent, \item[(iii)]
$r$ is a zerodivisor, \item[(iv)] $r$ is not invertible,
\item[(v)] $r$ belongs to every prime ideal in $R$.
\end{itemize}
\end{proposition}

\ \\
\subsection{$X(R) \neq \emptyset$}
In this case $I$ is not the only (maximal) prime ideal in $R$ and
$R$ is not a local ring.\\ For every prime $p$ and every signature
ideal $\mathcal{P} \in X(R)$, $\mathcal{P}+pR$ is
another maximal ideal.\\

For a given signature ideal $\mathcal{P}$ the intersection,
ranging over all primes p, of $\mathcal{P}+{pR}$ is just
$\mathcal{P}$. This implies that
\[{\rm Jac}(R) = \bigcap_{\mathcal{P}\in X(R),p \text{ prime}}\mathcal{P}+{pR}=\bigcap_{\mathcal{P}\in X(R)}\mathcal{P}={\rm Nil}(R).\]

First remark that $R_t\subseteq {\rm Nil}(R) = {\rm Jac}(R)$.
Suppose that $r \in R_t$ and $r \notin \mathcal{P}+{pR}$ for some
$\mathcal{P} \in X(R),\ p\ $prime. Consider the ideal
$J=(\mathcal{P}+{pR} +rR)$. Since $\mathcal{P}+{pR}$ is maximal,
$J=R$. So, $1 = s + rt$, for some $s \in \mathcal{P}+{pR}$, $t \in
R$. Let $k \in \mathbb{N}$ such that $kr = 0$, then we have $k =
ks \in \mathcal{P}+{pR}$. Since $\mathcal{P}+{pR}$ has finite
index $p$, this implies $p|k$, and so
$r \in \mathcal{P}+{pR}$, a contradiction.\\

Every nilpotent element is even ($l(r^k)\equiv l(r)^k \mod 2$).
Suppose that $r \in {\rm Nil}(R)$ and $l(r)=n$, $\ n$ even. Then
$p_n(r)=0$ i.e.
\[ r^{n+1}+c_{n-1}r^{n-1}+\ldots+c_3r^3+c_1r=0,\ \text{where}\ c_1\neq 0.\]
If $r^k=0$, $r^{k-1}\neq 0$ then multiplying the equation by
$r^{k-3}$ yields $c_1 r^{k-2} = 0$.  Multiplying the equation by
$c_1 r ^{k-5}$ then yields $c_1^2r^{k-4}=0.$  Repeating this
process will give $c_1^lr=0$ for some $l\in \mathbb{N}$. So $r \in R_t$.\\
This concludes that
\[ R_t = {\rm Nil}(R) = {\rm Jac}(R) = \bigcap_{\mathcal{P}\in X(R)}\mathcal{P}.\]

The unit group $\ R^\times = \{ r\in R\ |\ \sigma(r)=\pm 1,\
\text{ for all signatures } \sigma \text{ of }R\}$. This can be seen as follows:\\
If $r \in R^\times$ then there exists an $s\in R$ such that $rs=1$.
This implies that $\sigma(r)\sigma(s)=1$ and since $\sigma(r) \in
\mathbb{Z}$ we have $\sigma(r)=\pm 1$. On the other hand, suppose
that $\sigma(r)=\pm 1$ for all signatures $\sigma$ of $R$. Then we
have $\sigma(r^2 - 1) = 0$ for all $\sigma.$ So, $r^2 - 1 \in
\bigcap_{\mathcal{P}} \mathcal{P}={\rm Nil}(R)$. There exists a $k
\in \mathbb{N}$ such that $(r^2-1)^k=0$.  Now $(r^2-1)$ nilpotent
implies $r^2$ invertible, i.e. $r \in R^\times$.

Recall that for a commutative ring $R$, the set of zerodivisors,
${\rm Zd}(R)$, is the union of a certain set of prime ideals in $R$.
We will first show that $\bigcup_{\mathcal{P}\in
X(R)}\mathcal{P}\subseteq {\rm Zd}(R)$. Suppose that
$r\in\mathcal{P}\in X(R)$. Then $p_n(r)=0$ for some $n$ even. This
implies that $r[(r-n)\ldots(r-2)(r+2)\ldots(r+n)]=0$. Since
$\sigma((r-n)\ldots(r-2)(r+2)\ldots(r+n))\neq 0$ we have $r \in
{\rm Zd}(R)$.\\
Denote by $R_{t,p}$ the subset of $R_t$ consisting of the
$p$-torsion elements of $R$.  If $R_{t,p}\neq \{\,0\}$, then $p \in
\mathcal{P}+{pR}$ is a zerodivisor, implying that
$\mathcal{P}+{pR}\subseteq {\rm Zd}(R)$. Since all the signature
ideals $\mathcal{P} \subseteq \mathcal{P}+{pR}$, we have ${\rm
Zd}(R) = \bigcup_{\mathcal{P}\in X(R)}\mathcal{P}$ when $R_t=
\{\,0\}$ and ${\rm Zd}(R) = \bigcup_{\mathcal{P}\in X(R),p \text{
prime}}\mathcal{P}+{pR}$ for all $p$ prime such that $R_{t,p}\neq
\{\,0\}$. Since $\mathcal{P} \subseteq I = \mathcal{P}_{2}$ for all
$\mathcal{P} \in X(R)$, we have ${\rm Zd}(R) \subseteq I$ and the
same arguments hold as in the case $X(R) = \emptyset$ to prove that
$R$ is a 'connected' ring, i.e. $0$ and $1$ are the only
idempotents.\\
\ \\
These results can be summarized in the following propositions:
\begin{proposition}
Let $R$ be an admissible AP ring with generating polynomial $q(X)=X^2-1$ and assume that $X(R)\neq\emptyset$.\\
Then for an element $r \in R$, the following conditions are
equivalent :\\
\begin{itemize}
\item[(i)] $r$ is a torsion element , \item[(ii)] $r$ is
nilpotent, \item[(iii)] $r$ belongs to every prime ideal in $R$,
\item[(iv)] $r$ belongs to every signature ideal in $R$.
\end{itemize}
The equivalence $(i)\Leftrightarrow(iv)$  is in fact Pfister's
local-global principle.
\end{proposition}
The set of zerodivisors is completely described by the following
\begin{proposition}
Let $R$ be an admissible AP ring with generating polynomial $q(X)=X^2-1$ and assume that $X(R)\neq\emptyset$.\\
If $R$ is torsion-free then ${\rm Zd}(R)=\bigcup_{\mathcal{P} \in
X(R)}\mathcal{P}$, the union of all signature ideals in $R$.\\
Otherwise, ${\rm Zd}(R)=\bigcup_{\mathcal{P}\in
X(R),p}(\mathcal{P}+pR)$,  for $p$ prime such that $R$ has non-zero $p$-torsion.\\
\end{proposition}
\section{Constructing annihilating polynomials}
We will construct annihilating polynomials for AP rings for
different choices of the generating polynomial $q(X)$.
\subsection{$q(X)=X^2-1$}
This is the well-known case described in \cite{LEW90}. The roots
of the generating polynomial are $-1$ and $1$.  The possible
values for the sum of $n$ elements out $\{-1,1\}$ lie in
$\{-n,-n+2,\ldots,n-2,n\}$. The annihilating polynomial for an
element of length $n$ is thus the $n$-th Lewis polynomial
\[p_n(X) = (X-n)(X-(n-2))\ldots(X+(n-2))(X+n).\]
\subsection{$q(X)=X^4-1$}
Write $R_1=\{-1,1,-i,i\}$ for the set of roots of $q(x)=x^{4}-1$ in
$\mathbb{C}$. Denote by $R_j$ the subset of the complex numbers
$\mathbb{C}$ consisting of sums of $j$ elements of $R_1$. Since
$0\in R_2$ we have $R_{n}\subset R_{n+2}$ for all $n \in
\mathbb{N}^\ast$. Consider $D_n=R_{n}\setminus R_{n-2}$ for $n>2$.
Put $D_1=R_1$ and $D_2=R_2$.\\ Now define the monic integer
polynomial $t_n(X) \in \mathbb{Z}[X]$ as follows:
\[t_n(x)=(x^4-n^4)\prod_{\substack{a+b=n\\a,b \in \mathbb{N}^\ast}}(x^4-2(a^2-b^2)x^2+(a^2+b^2)^2)\]
\begin{align*}
t_1(x)=&\ x^4-1\\
t_2(x)=&\ (x^4-16)(x^4+4)\\
t_3(x)=&\ (x^4-6x^2+25)(x^4+6x^2+25)\\
t_4(x)=&\ (X^4-256)(x^4-16X^2+100)(x^4+64)(x^4+16X^2+100)\\
\end{align*}
\begin{lemma}
The polynomial $t_n(x)$ has the property that $t_n(z)=0$ for all
$z\in D_n$.
\end{lemma}
\begin{proof}
For $n>1$, this follows from the observation that \[D_n=\{a+bi \in
\mathbb{Z}[i]\,:\,|a|+|b|=n\}.\]
\end{proof}
Now we are able to construct the annihilating polynomial $p_n(x)$.
Since $R_{n}\subset R_{n+2}$ and $0\in R_{2n}$ for all $n \in
\mathbb{N}^\ast$ it follows that:
\[\text{For $n$ even }p_n(x)=t_n(X)t_{n-2}(X)\ldots t_2(X)X.\]
\[\text{For $n$ odd }p_n(x)=t_n(X)t_{n-2}(X)\ldots t_1(X).\]
A few examples:
\begin{align*}
p_1(x)=&\ x^4-1\\
p_2(x)=&\ x(x^4-16)(x^4+4)\\
p_3(x)=&\ (x^4-1)(x^4-6x^2+25)(x^4+6x^2+25)\\
p_4(x)=&\ x(x^4-16)(x^4+4)(X^4-256)(x^4-16X^2+100)(x^4+64)(x^4+16X^2+100)\\
\end{align*}
Elements in the Witt ring of level $2$ of dimension $n$ will be
annihilated by $p_n(X)$.
\subsection{$q(X)=X^{2^k}-1$}

Let $R$ be an AP ring with generating polynomial $q(X)=X^{2^k}-1$.
The set of roots of $q(X)$ is generated by a primitive $2^k$-th
root of unity $\zeta$. In general it is very difficult to find an
explicit expression for the annihilating polynomial $p_n(X)$ but
we will obtain some properties.  We prove the following
\begin{lemma}\label{odd constant term}
Let $R$ be an AP ring with generating polynomial $q(X)=X^{2^k}-1$
and let $p_n(X)$ be the annihilating polynomial. Then $p_n(0)$ is
odd whenever $n$ is odd.
\end{lemma}
\begin{proof}
Let $R_1=\{\zeta,\zeta^2,\ldots,\zeta^{2^k}\}$ be the set of roots
of $q(X)=X^{2^k}-1$ in $\mathbb{C}$ and let $R_i$ be the subset of
the complex numbers $\mathbb{C}$ consisting of sums of $i$
elements of $R_1$.  Note that \[R_n = \{\sum_{i=1}^{2^k}a_i\zeta^i
\quad|\quad a_i \in
\mathbb{N}^\ast\quad\text{and}\quad\sum_{i=1}^{2^k}a_i=n\}.\] Put
$P:=\prod_{\sigma\in R_n}\sigma$.\\
Since $\zeta^j=-(1-\zeta)(1+\zeta+\ldots+\zeta^{j-1})+1$ it follows
that \[ P = n^l + (1-\zeta)m(\zeta)\] for some complex function $m$
and $l=\#R_n$. Then \[ P^{2^k} = n^{l\cdot2^k} + 2m'(\zeta)\] for
some rational valued function $m'$ (since $P \in \mathbb{Z}$). It
follows that $P$ is odd whenever $n$ is odd. The result now follows,
since $p_n(0) \mid P$.
\end{proof}
The following lemma will provide us an upper bound on the degree
of the annihilating polynomials.
\begin{lemma}
Let $R$ be an AP ring with generating polynomial $q(X)=X^{2^k}-1$
and let $p_n(X)$ be the annihilating polynomial. Then ${\rm deg}\
p_n(X)\leq 2^{n-1}(2^k-1)+1$.
\end{lemma}
\begin{proof}
To determine the degree of the annihilating polynomial $p_n(X)$ we
have to calculate the number of different elements in the set
\[R_n = \{\sum_{i=1}^{2^k}a_i\zeta^i \quad|\quad a_i \in
\mathbb{N}^\ast\quad\text{and}\quad\sum_{i=1}^{2^k}a_i=n\}.\] We
write
\begin{align*}
\sum_{i=1}^{2^k}a_i\zeta^i =&
\sum_{i=1}^{2^{k-1}}(a_i-a'_i)\zeta^i\\
&\qquad\qquad\qquad\text{where }a'_i=a_{2^{k-1}+i}\\
=&\sum_{i=1}^{2^{k-1}}b_i\zeta^i\qquad\text{(1)}\\
&\qquad\qquad\qquad\text{where }b_i>0\\
&+\sum_{i=1}^{2^{k-1}}b'_i\zeta^i\qquad\text{(2)}\\
&\qquad\qquad\qquad\text{where }b'_i<0\\
&+0.
\end{align*}
The number of elements in (1) and (2) are for symmetry reasons the
same and it follows that $\#R_n = 2 \times \#R_{n-1} - 1$.\\
Since $\#R_1=2^k$, we obtain that ${\rm deg}\ p_n(X) \leq \#R_n =
2^{n-1}(2^k-1)+1.$
\end{proof}
An example of an AP ring with generating polynomial $X^{2^k}-1$ is
a Witt rings of level $k$ as defined in \cite{Kleinstein
Rosenberg}.
\subsection{$q(X)=X^2-2^kX$}
Write $R_1=\{0,2^k\}$ for the set of roots of $q(x)=x^{2}-2^kX$ in
$\mathbb{C}$. Denote by $R_n$ the subset of the complex numbers
$\mathbb{C}$ consisting of sums of $n$ elements of $R_1$. So
$R_n=\{0, 2^k,2\cdot2^k,\ldots,n\cdot2^k\}$ and $p_n(X)$  is given
by
\[p_n(X)=X(X-2^k)(X-2\cdot2^k)\ldots(X-n\cdot2^k).\]

An example of an AP ring with generating polynomial $X^2-{2^k}X$
is for example the subring of the Witt ring additively generated
by all $k$-fold Pfister forms of this Witt ring.  In this setting,
$p_n(X)$ will annihilate all sums of $n$ $k$-fold Pfister forms.

\subsection{Table of Marks}
The {\emph table of marks} of a finite group $G$, as introduced by
Burnside in his classic theory of groups of finite order \cite{BUR},
is a matrix whose rows and columns are labeled by the conjugacy
classes of subgroups of $G$ and where for two subgroups $H_1$ and
$H_2$ the $(H_1,H_2)$-entry is the number of fixed points of $H_2$
in the transitive action of $G$ on the cosets of $H_1$ in $G$. This
table of marks is sometimes called a \emph{ Burnside matrix}.  If
$n_1,n_2,\ldots,n_k \in \mathbb{N}$ are the different entries in the
table of marks of a group $G$, then, by example $2.2.(viii)$, the
polynomial $q(X)=\prod_{i=1}^k(X-n_i)$ is a generating polynomial
for the Burnside ring $\Omega(G)$.
\begin{example}
The table of marks of the alternating group $A_5$ is given by
\[
\begin{array}{lccccccccc}
&e&C_2&C_3&V_4&C_5&S_3&D_{10}&A_4&A_5\\
&&&&&&&&&\\
e&60&0&0&0&0&0&0&0&0\\
C_2&30&2&0&0&0&0&0&0&0\\
C_3&20&0&2&0&0&0&0&0&0\\
V_4&15&3&0&3&0&0&0&0&0\\
C_5&12&0&0&0&2&0&0&0&0\\
S_3&10&2&1&0&0&1&0&0&0\\
D_{10}&6&2&0&0&1&0&1&0&0\\
A_4&5&1&2&1&0&0&0&1&0\\
A_5&1&1&1&1&1&1&1&1&1\end{array}
\]
The generating polynomial for the AP ring $\Omega(A_5)$ is thus
given by
$q(X)=(X-60)(X-30)(X-20)(X-15)(X-12)(X-10)(X-6)(X-5)(X-3)(X-2)(X-1)X$.
\end{example}
%%%%%%%%%%%%%%%%%%%%%%%%%%%%%%%%%%%%%%%%%%%%%%%%%%%%%%%%%%%%%%%%%%%%%%%%%%%%%%%%%%%%%%%%%%%%%%%%%%%%%%%%%%%%%%%%%%
%                                                                                                                %
%  THE BIBLIOGRAPHY                                                                                              %
%                                                                                                                %
%%%%%%%%%%%%%%%%%%%%%%%%%%%%%%%%%%%%%%%%%%%%%%%%%%%%%%%%%%%%%%%%%%%%%%%%%%%%%%%%%%%%%%%%%%%%%%%%%%%%%%%%%%%%%%%%%%

\noindent David W. Lewis \\
School of Mathematical Sciences \\
University College Dublin \\
Belfield, Dublin \\
Ireland \\
e-mail: {\tt david.lewis@ucd.ie}\\
\ \\
\noindent Stefan De Wannemacker \\
Faculty of Applied Economic Sciences \\
University of Antwerp \\
Prinsstraat 13, 2000 Antwerp \\
Belgium \\
e-mail: {\tt stefan.dewannemacker@ua.ac.be}
\end{document}